\input amstex.tex

\input amsppt.sty

\TagsAsMath

\NoRunningHeads

\magnification=1200

\hsize=5.0in\vsize=7.0in

\hoffset=0.2in\voffset=0cm

\nonstopmode

\document

\def\pr{^\prime }

\document
\topmatter

\title{On dispersion   for Klein Gordon equation with
periodic potential in 1D}
\endtitle

\author
Scipio Cuccagna
\endauthor

\address
DISMI University of Modena and Reggio Emilia, Reggio Emilia 42100
Italy
\endaddress

\email cuccagna.scipio\@unimore.it\endemail

\abstract By exploiting estimates on Bloch functions obtained   in a
previous paper, we prove   decay estimates for
  Klein Gordon equations with a time independent potential
periodic in space in 1D and with  generic  mass.

\endabstract

\endtopmatter

\head \S 1 Introduction  \endhead

We consider    Schr\"odinger operators of the form $H  =H_0 +P(x)$
with, $H_0= -\frac{d^2}{dx^2}$, $P(x)$ a smooth nonconstant real
valued periodic function, $P(x+1)\equiv P(x)$,  with spectrum
$\Sigma (H )=\cup _{n\ge 0} \Sigma _n$, formed by bands $\Sigma
_n=[A_n^+, A_{n+1}^-]$ with $A_n^+< A_{n+1}^-\le A_{n+1}^+$ for any
$n\in \Bbb N \cup \{ 0 \}$. We normalize $H $ so that $A_0^+=0$. We
then show:

\proclaim{Theorem 1.1} Under the above hypotheses consider  for
  $\mu >0$  the solutions of the following Cauchy problem
for the Klein Gordon equation
$$u_{tt}+H u+\mu u=0\, , \quad u(0,x)\equiv 0 \, , \quad
u_t(0,x)=g(x). \tag 1.1$$ Then there exists a bounded discrete set $
\Bbb D\subset (0,+\infty )$   such that for any $\mu \in (0,+\infty
) \backslash\Bbb D$
 there is  a   $C_\mu >0$
such that the  following dispersive estimate holds:
$$\| u(t,\cdot )\| _{L^\infty (\Bbb R)}\le C_\mu  \langle t \rangle ^{-\frac{1}{3}}\| g( \cdot )\| _{W^{ 1,1}(\Bbb R)}
.\tag 1.2$$
\endproclaim
Maybe $ \Bbb D$ is empty. The exact condition defining $ \Bbb D$ is
given in Lemma 3.1 below. The proof is based on results in \cite{C}
where proofs are explicitly done only for the generic case when all
the spectral gaps are nonempty. Since the generic case contains all
the crucial difficulties, there is no problem at extending the
results in \cite{C} to the non generic case, and we will assume this
as a fact (and if this is unconvincing the reader can assume that
$\ell _n=n$ below). To illustrate Theorem 1.1 consider $P(x)=2\kappa
^2 \text{sn} ^2 (x,\kappa )$ for $\kappa \in (0,1)$, with $\text{sn}
(x,\kappa )$ the Jacobian elliptic function. Then $\Sigma (H
)=[\kappa ^2,1]\cup [ 1+\kappa ^2, +\infty ) $ and  by Theorem 1.1
for generic $\mu
>-\kappa ^2$ we get (1.2). Notice that for $\mu=0$   this example  resembles the
flat Klein Gordon rather than the flat wave equation, because we
have  $A_0^+ =\kappa ^2>0$.   For $H=H_0$ the equation $u_{tt}+H u -
|u|^p=0$ for any $p>1$ is not globally well posed for small initial
data in $C^\infty _0(\Bbb R)$ while if $p\gg 1$ this is not the case
for $H$ with $P(x)= 2\kappa ^2 \text{sn} ^2 (x,\kappa )$ or
$P(x)=\sin ^2(x)$. In the latter case all the gaps are non empty.
The proof in this paper mixes results from \cite{C} with a specific
computation in Marshall {\it et al.} \cite{MSW}, specifically Lemma
5 therein.

  \head \S 2 Reformulation, spectrum,  band and Bloch functions \endhead

We will prove:

\proclaim{Theorem  2.1} Let $H $ be as in Theorem 1.1, that is with
a smooth periodic potential, and such that $A_0^+=0.$ Then, there is
a set $\Bbb D$   like in Theorem  1.1 such that for any $\mu \in
(0,+\infty ) \backslash \Bbb D$
 there is  a   $C_\mu >0$ such that the  following
dispersive estimate holds:
$$\left  \|  {\sin (t\sqrt{ H  +\mu } )}{(H  +\mu )^{-\frac{3}{4}}} \colon L^1(\Bbb R)
\to L^\infty (\Bbb R)\right  \| \le C _\mu \langle t \rangle
^{-\frac 13}
  .\tag 2.1$$
\endproclaim

\noindent The $u$ in  (1.1) is, for $1/4>\varepsilon >0$,  $ G=
 (H  +\mu )^{ \frac{1}{4}}  (H_0 +1
)^{-\frac{1}{2}+\varepsilon} $ and $h=(H_0 +1 )^{
\frac{1}{2}-\varepsilon}g $, $u(t)=$  $$ \frac{\sin (t\sqrt{ H +\mu
} ) }{ (H  +\mu )^{ \frac{3}{4}}} G h \Rightarrow  \| u(t)\|
_{\infty}\le \left \| \frac{\sin (t\sqrt{ H +\mu } )}{(H  +\mu
)^{\frac{3}{4}}} \right \| _{L^1\to L^  \infty} \left \| G \right \|
_{ L^1\to L^  1} \| h\| _{1} . $$ We have $\| h\| _{1} \le C \| g\|
_{W^{1,1}}$, $ \left \|G  \right \| _{L^1\to L^1}\le C(\mu )$, with
$C(0)=O(1)$, so (2.1) implies (1.2). (2.1) is a consequence of the
following estimate:

\proclaim{Proposition  2.2} There is a set $ \Bbb D$   like in
Theorems 1.1-2 such that for any  $\mu \in (0,+\infty ) \backslash
\Bbb D$
 there is  a   $C_\mu >0$ such that  the  following estimate holds
for any $(t,x,y)$:
$$\left  | \left ({\sin (t\sqrt{ H  +\mu } )}{(H  +\mu )^{-\frac{3}{4}}} \right )(x,y)\right  |
\le C _\mu \langle t \rangle ^{-\frac 13}
  .\tag 2.2$$
\endproclaim
We will prove Proposition 2.2 in the case when the spectrum as
$\Sigma (H )$ is formed by infinitely many bands, the finitely many
bands case being easier. To prove (2.2) we    express   the integral
kernel in (2.2) in terms of   Bloch functions, see below. We express
 $\Sigma (H )=\cup _{n=0}^\infty \Sigma _n$,
with $\Sigma _n=[A_n^+, A_{n+1}^-]$ with $A_n^+< A_{n+1}^-\le
A_{n+1}^+$ for any $n\in \Bbb N \cup \{ 0 \}$, with $A_0^+=0.$ Set
for $n\ge 0$, $a_{n}^{\pm }=\sqrt{A_{n}^{\pm }}$ and $a_{-n}^{\pm
}=-a_{n}^{\pm }$. For $n\ge 0$ set $\sigma _n=[a_n^+, a_{n+1}^-] $
and $\sigma _{-n}=-\sigma _n$.  Set $ \sigma = \cup
_{n=-\infty}^\infty \sigma _n$, with each two   intervals $\sigma
_n$  and $\sigma _{n+1}$ separated by a non empty  gap $g_n$. For
$|g_n|$ the length of the gap $g_n$ we have   the following
classical result, see \cite{E} ch. 4:

\proclaim{Theorem 2.3} Let $P(x)$ be smooth. Set $\sigma _n=[a_n^+,
a_{n+1}^-]$ and $g_n=]a_n^-, a_{n }^+[$. Then $\exists$ a strictly
increasing
  sequence $\{ \ell _n \in \Bbb Z \} _{n\in \Bbb Z}$ and a
fixed constant $C $ such that
$$|a_n^- -\ell _n\pi |+ |a_n^+ -\ell _n\pi |\le C  \langle \ell _n \rangle ^{-1}.
$$
$\forall$ $N$ $\exists$ a fixed constant $C_N $   such that
$|g_n|\le C _N \langle \ell _n \rangle ^{-N}  $ $\forall $ $n$.
\endproclaim

We    review band   and Bloch functions.
 $\forall$ $w
\in {\Bbb C_+}$ (the open upper half plane)  $\exists$ a unique   $
k\in {\Bbb C_+}$ such that there are two solutions of $(H-w^2)u=0$
of the form  $\tilde \phi _\pm (x,w)=e^{\pm ikx} \xi _\pm (x,w)$
with  $\xi _\pm (x+1,w)\equiv \xi _\pm (x,w)$ and with $\tilde \phi
_\pm (0,w)=1$.  The correspondence between $w$ and the
"quasimomentum" $k$ is a conformal map  between ${\Bbb C_+}$  and a
"comb" $K ={\Bbb C_+}-\cup _{n\neq 0}  [  \ell _n\pi ,  \ell _n\pi
+ih _n ]$  with $\ell _n$  satisfying  the conclusions of Theorem
2.3,   with $|g_n|\le 2 h_n \le C |g_n|$ for a fixed $C$. For
generic potentials, $\ell _n\equiv n$. The map $k(w)$ extends into a
continuous map in $\overline{\Bbb C_+}$ with $k(\sigma _{ _n})=[
\ell _{n}\pi , \ell _{n+1}\pi ]$, with $k(w)$ a one to one and onto
map between $\sigma _n$ and $[ \ell _{n}\pi , \ell _{n+1}\pi ]$, and
with
 $k(g _n)=] \ell _{n}\pi ,
\ell _{n}\pi +ih_n ]$. $k(w)$ extends into a conformal map from
$\Bbb C -\cup _{n\neq 0}\overline{g _n}$ into $\Cal K=\Bbb C-\cup _n
\gamma _n $  with $\gamma _n= [\ell _{n}\pi -ih_n, \ell _{n}\pi
+ih_n]$. Next set $N^2(w)= \int_0^{1} \tilde \phi _+ (x,w)\tilde
\phi _- (x,w) dx.$ We have $N^2(w)= \int_0^{1}   \big |\tilde  \phi
_\pm (x,w) \big |^2 dx>0$ for $w\in \sigma $, $N^2(w)\neq 0$ for any
$w\in \Bbb C \backslash  \cup _{n\neq 0}\overline{g _n}$. We set $ m
_+ ^0(x,w)m _- ^0 (y,w)=  {\xi _+(x,w)\xi _-(y,w)}{N ^{-2 } (w)}
 $ and   $m _\pm ^0(x,w)=\xi _\pm ^0(x,w)/N(w)$  with
$N(w)>0$ for $w\in \sigma$. We express $w=w(k)$ for $k\in \Cal K$
and with an abuse of notation we write $\phi _\pm (x,k)$ for $\phi
_\pm (x,w(k))$ and $m^0 _\pm (x,k)$ for $m^0 _\pm (x,w(k))$.   We
call $\phi _\pm (x,k)= e^{\pm ikx}\phi _\pm (x,k)$ Bloch functions.
 In \cite{C}  we
had to work with $w$ complex, but here  we focus only on $w\in\sigma
$.  The band function is $E(k)=w^2(k)$. Now we have the following
well known fact:

\proclaim{Theorem 2.4} Set $\hat f (k)=\int _\Bbb R  {\phi  _+
(y,k)} f(y)dy$ for any $k\in \Bbb R\backslash \pi \Bbb Z$. Then:

$$\align & \int  _\Bbb R| f(y)|^2dy=  \int  _{\Bbb R}|\hat f (k)|^2dk
, \quad   f(x)=   \int _{\Bbb R}  {\phi _- (x,k)} \hat f(k) dk,\quad
\widehat{H f  }(k)= E (k)\hat f(k)  .   \endalign
$$
In particular we have $$ \frac{\sin (t\sqrt{ H  +\mu } )}{( H  +\mu
) ^{\frac{3}{4}}} (x,y)=
  \int   _{\Bbb R}
e^{ -i(x-y)k  } \frac{\sin (t\sqrt{ E  (k)  +\mu } )}{( E (k) +\mu )
^{\frac{3}{4}}}  {m _- ^0(x,k)}{ m  _+ ^0(y,k)} dk .\tag 2.4$$

\endproclaim
We will show that the generalized integral (2.4) is a function which
satisfies (2.2).

  \head \S 3 Estimates on band and Bloch functions \endhead

  We set  $\dot f=
\frac{df}{dk}$,  $f'= \frac{df}{dw}$ and $\eta (k)= \sqrt{E (k) +\mu
}$. We compute

$$\aligned & \dot \eta  = \frac{\dot E}{2(E+\mu )^{\frac{1}{2}}}  \, , \quad  \ddot \eta  = \frac{\ddot E}{2(E+\mu
)^{\frac{1}{2}}} -\frac{\dot E^2}{4(E+\mu )^{\frac{3}{2}}};
  \\& \dddot \eta  = \frac{\dddot
 E}{2(E+\mu )^{\frac{1}{2}}} -\frac{3\dot E \ddot
 E}{4(E+\mu )^{\frac{3}{2}}} +\frac{3\dot E ^3}{8(E+\mu )^{\frac{5}{2}}} =\frac{\dddot
 E}{2(E+\mu )^{\frac{1}{2}}}-\frac{3}{2} \dot E \ddot \eta  .
\endaligned \tag 3.1$$

 \proclaim{Lemma 3.1} $\exists$   $
\Bbb D\subset (0,+\infty )$, bounded and discrete,   such that
$\forall$ $\mu \in (0,+\infty )\backslash \Bbb D$   the system
$\ddot \eta (k)=\dddot \eta (k)=0$, or equivalently  (3.2) below,
has no solutions in $\Bbb R$:
$$  \ddot E=  \frac{\dot
E^2}{2(E+\mu ) }  \, , \quad  \dddot
 E= 0.\tag 3.2$$

\endproclaim
For the case $A^+_0=0$ and  $\mu =0$ see  Korotyaev \cite{K1}.

{\it Proof of Lemma 3.1.} We start by focusing on {\it low
energies}.   $|E|\le E_0$ implies $|k|\le k_0$ for a fixed
$k_0=k_0(E_0)$. By \cite{K2} we have the following two facts:

\proclaim{Lemma 3.2} (a) On each band,  $\dot E=0$ holds exactly at
the extremes of the band.

  \noindent (b) On each   band, there is exactly
one solution of  $\ddot E=0$, contained in the interior of each
band.
\endproclaim
Here recall we are assuming the bands to be bounded.  By (a), for
$|k|\le k_0$ (3.2) cannot hold near the extremes of the bands. So
there is a fixed $c>0$, such that, if $k$ is a root of (3.2), then
$k$ is in the set, which we denote by $J$, formed by the   $k$ whose
distance from the nearest edge is at least $c$.
  $\dddot E=0$ has finitely many solutions in $J$. Indeed, $\dddot
  E\not \equiv 0$, is holomorphic in $\Cal K$ and
  $\overline{J}\subset \Cal K$.
 So except for at most finitely $\mu $'s with $\mu >0$,
(3.2) has no solutions for $|k|\le k_0$.

 Next we consider Lemma 3.1 in the   {\it high energy}
 case.    Recall  $E=w^2$ and  assume that (3.2)
is satisfied at some   value   $w_0$. Since $E$ is even we can
assume $w_0>0$, in particular $w_0\in [a^+_n, a^-_{n+1}]$. We have
$a^+_n=w((\ell _n   \pi )^+ )$, with $\ell _n\in \Bbb N$. We have
the following facts:

\proclaim{Lemma 3.3} (1) There is a fixed $C>0$ such that for $
a_{n+1 }^--C
 \ell  _{n+1}^{\frac{1}{3}}|g_{n+1 }|^{\frac{2}{3}} <w\le a^-_{n+1}$
 we have $\ddot E<0$.

 \noindent (2) For any given $C_1\gg 1$ there are
 $n_0$ and $c_0>0$ such that for $n\ge n_0$ and for
$$a_{n+1 }^--C_1
\ell  _{n+1}^{\frac{1}{3}}|g_{n+1 }|^{\frac{2}{3}} \le w
 \le  a_{n+1 }^--C
 \ell  _{n+1}^{\frac{1}{3}}|g_{n+1 }|^{\frac{2}{3}} $$
 we have $|\dddot E|\ge c_0
 \ell  _{n+1}^{-1 } |g_{n+1}|^{-\frac{2}{3}} \gg 1 $.

 \noindent (3)   If $a^-_{n+1}< \infty $ then there exists exactly one
   point $w_1$,  $w_1\in(a^+_n,
 a^-_{n+1})$,  with $k''(w_1)=0$. For   $w\in [a^+_n,w_1)$ we have $k''(w)<0$. Furthermore, there are positive constants $C_0$, $C_1$, $C_2$,
   $\alpha $ such
  that    for any $ n\ge n_0$,   for any $w> a^-_{n+1} -\alpha $
 we have

$$\align &   \frac{C_1}{\langle w
\rangle ^{ 3}}-\frac {\varphi  (w)}4\ge  - k^{\prime \prime }(w) \ge
\frac{C_2}{\langle w \rangle ^{ 3}}-\frac {C_0\varphi  (w)}4, \text{
$\varphi  (w):= \frac{(a_{n+1}^+-a_{n+1}^-)^2}{|w-a_{n+1}^-|^{\frac
32}|w-a_{n+1}^+|^{\frac 32}}$;}
\endalign$$
\noindent (4) For $ a_{n  }^-+ |g_{n }|^{\frac{1}{4}} \le w
 \le  a_{n+1 }^-- |g_{n+1 }|^{\frac{3}{5}} $
we have  $|\dot w -1+ \frac{Q_0}{k^2}|\le C |k |^{-3} $ for a fixed
$C$, with $Q_0:=\frac{1}{2}\int _0^1P(x) dx$. Furthermore, $\dot E=
2k+O(k^{-2})$ and $E=k^2+2Q_0+ O(k^{-2})$ with in either case
$|O(k^{-2})|\le C k^{-2}$ for a fixed $C$.

\endproclaim
(1) is a consequence of Lemma 4.2 \cite{C}, (2) is   Lemma 4.3
\cite{C}. In (3) the information on the sign of $k''(w)$ is in
\cite{K1} and the inequality is in Lemma 7.3 \cite{C}. In (4) the
inequalities for $\dot w$ and $\dot E$ are proved in Lemma 7.1
\cite{C}, the inequality for $ E$ is proved in Lemma 5.4 \cite{C}.

\noindent We return to Lemma 3.1.
   $ \ddot E>0$  at $w_0$ by (3.2). By
(1) Lemma 3.3,
    $    w_0<a_{n+1 }^--C
 \ell  _{n+1}^{\frac{1}{3}}|g_{n+1 }|^{\frac{2}{3}} $. By $\dddot E=0$ at $w_0$   and by (2) Lemma 3.3 then $ w_0
 \le  a_{n+1 }^--C_1
 \ell  _{n+1}^{\frac{1}{3}}|g_{n+1 }|^{\frac{2}{3}} $ for some $C_1\gg
 1$. By (3) Lemma 3.3 in this region $k''<0$.
Hence we have the inequality
 $  \ddot E=2(\dot w)^2-2w(\dot w)^3k''\ge  2(\dot w)^2.$
Suppose   $a^{+}_{n}+C_1 \ell  _{n } ^{3} |g_n| \le w_0$. By (4)
Lemma 3.3 we have $ \dot w=1- \frac{Q_0}{k^2}+O(k ^{-3})$. So
$$\ddot E \ge 2 ( 1- 2\frac{Q_0}{k^2}) +O(k^{-3}) .\tag 1$$ On the
other hand   for $w_0=w(k )$ by (4) Lemma 3.3 we have
$$\ddot E=  \frac{\dot
E^2}{2(E+\mu  ) }= 2 \frac{( k + O(k ^{-2}))^2 }{   k^2 + 2Q_0+ \mu
+ O(k^ {-2} )} = 2- \frac{4Q_0+\mu }{k^2 } +O(k^{-3} ) $$ The last
  formula  is incompatible with   (1)
for $\mu \ge \mu _0>0$ with $\mu _0$ fixed and for $|k| \gg 1/\mu $.
Hence   at large energies   and for $a^{+}_{n}+C_1 \ell _{n } ^{3}
|g_n| \le w $, for some fixed $C_1>0$, there are no solutions of
(3.2). Let  $a^{+}_{n}\le w\le a^{+}_{n}+C_1 \ell  _{n } ^{3}
|g_n|.$ The following lemma, see Lemma 4.3 \cite{C},  shows $w_0\not
\in \left [ a^+_{n }+c|g_n|,a^+_{n }+|g_n|^{\frac{3}{5}} \right ]
 $  for $c \gg 1$ fixed:

 \proclaim{Lemma 3.4} For $a^+_{n }+c|g_n|\le w\le a^+_{n }+|g_n|^{\frac{3}{5}}
 $ for $c\gg 1$ a fixed large constant,
$|\dddot E|$ is very large. \endproclaim

 Finally $w_0\not \in  \left
[a^+_{n }, a^+_{n }+c|g_n|  \right ]
 $ because by the following lemma, see Lemmas 7.1
and 7.4 \cite{C},  and by Lemma 2.3 for $|u -a^+ _{n}|\lesssim |g
_{n}| $ then $ \ddot E=  \frac{\dot E^2}{2(E+\mu ) }$  cannot hold:

\proclaim{Lemma 3.5} For $|u -a^+ _{n}|\le c |g _{n}| $  for $c>0$
fixed there are $n_0$,     $C_1>0$ and $C_2>0$  such that  for any
$n\ge n_0$ we have $|\ddot E|\ge C_1\ell _n|g_n|^{-1}$ and $|\dot
E|\le C_2 \ell _n \sqrt{\frac{u-a^+_n}{|g_n|}}.$

\endproclaim
Finally for later use we state the following, see Lemma 7.1
\cite{C}:

\proclaim{Lemma  3.6} $\exists$ $C_1>C_2>0$ such that $\forall m$
and $\forall w\in \sigma _m=[a^+_m,a^-_{m+1}]$ we have for $A(w)=
\frac{ |g_m|^2 }{(
 w-a_m^+    )^{\frac{1}{2}}  ( w-a_m^++|g_m|)^{\frac{3}{2}}  }  +  \frac{ |g_{m+1}|^2
}{(
  a_{m+1}^- -w   )^{\frac{1}{2}}  (  a_{m+1}^--w+|g _{m+1}|)^{\frac{3}{2}}  }$
  $$ \aligned &  1+  C_2\left ( A(w) + \frac{1}{\langle w \rangle ^2} \right ) \ge
    k^\prime (w) \ge  1+C_1 A(w)    .\endaligned  $$
  Correspondingly for $k\in [\ell _m\pi , \ell  _{m+1}  \pi ]$ and for $\dot w= dw/dk$ we have
  $$   \align &\frac 1{1+C_2 \left ( A(w) + \frac{1}{\langle w \rangle ^2 } \right )}\le      \dot w\le
     \frac 1{1+C_1 A(w)}   . \endalign$$
\endproclaim

  \head \S 4 Decomposition of (2.3) and estimates on the single parts  \endhead

We decompose $  \frac{\sin (t\sqrt{ H +\mu } )}{(H+\mu )
^{\frac{3}{2}}} (x,y)=\sum _nK^n(t,x,y) $ with

$$K^n(t,x,y):=
  \int   _{[\ell _n \pi ,\ell _{n+1}\pi ]}e^{ -i(x-y)k  } \frac{\sin (t\eta (k) )}{\eta
^{\frac{3}{2}} (k)}  {m _- ^0(x,k)}{ m  _+ ^0(y,k)} dk. \tag 4.1$$ A
basic ingredient in the proof is the stationary phase theorem, see
p. 334 \cite{S}:

\proclaim{Lemma  4.1} Suppose $\phi (x)$ is real valued and smooth
in $[a,b]$ with $|\phi ^{(m)}(x) |\ge c_m>0$ in $]a,b[$ for $m\ge 1$
. For $m=1$ assume furthermore that $\phi ^\prime (x)$ is monotonic
in $]a,b[$. Then   we have for $C_m=  5 \cdot 2^{m-1}  -2 $:
$$ \left |\int _a^b   e^{i\mu \phi (x)}  \psi (x)  dx
\right | \le C_m  (c _m\mu )^{-\frac 1m} \left [ \min \{ |\psi (a)
|, |\psi (b) | \}+\int _a^b |\psi \pr (x) | dx \right ] . $$
\endproclaim

The following two lemmas are special cases of    Lemmas 4.4 \& 4.5
in \cite{C}:

\proclaim{Lemma 4.2} There  are fixed constants $C>0$, $C_3>0$,
  $\Gamma
>0$ and $c>0$ such that for all $x$, all $n$ we have :

{\item{ (1)}}   $  \forall \,   w \in [a^+_{n }
 +C_3\ell _{n }^5 |g_n| ,a^-_{n+1 }-C_3\ell _{n+1}^5 |g _{n+1}| ]$ we have
$  \big |  m ^0_+(x,k)  m ^0_-(y,k) -1 \big | \le  \frac{C}{\langle
k\rangle } ;
 $

{\item{ (2)}} $\exists$ fixed $C>0$ such that $forall$ $k\in \Bbb R
 $   we have
 $ \big |  m ^0_+(x,k)  m ^0_-(y,k)   \big | \le C .
 $

\endproclaim

  \proclaim{Lemma  4.3} There  are fixed
constants $C>0$ and $C_4>0$, with $C_4<C_2$,   such that for all
$x$, all $n$   we have:

{\item {(1)}} for all $a^{+}_{n}+ |g_n|^{\frac{1}{4}}\le k \le
\frac{a^{+}_{n}+a^{-}_{n+1}}{2} $, then
     $  | \partial _k (m_- ^0(x ,k) m_+
^0(y ,k))\big | \le \frac{C}{k |k-\pi \ell _{n }|}
 ;$

 {\item {(2)}} for all  $
 \frac{a^{+}_{n}+a^{-}_{n+1}}{2} \le k \le a^{-}_{n+1}-
 |g _{n+1}|^{\frac{3}{5}}$, then $  | \partial _k
(m_- ^0(x ,k)  m_+ ^0(y ,k)) \big | \le \frac{C}{k |k-\pi \ell _{n
+1}|}
 ;$

{\item {(3)}} for   $    w $ in the remaining part of $  [a^+_{n },
a^-_{n+1 }]$ we have for $m=n$ (resp. $m=n+1$)  near $a^+_{n }$
(resp. $a^-_{n+1 }$)
 $ \big | \partial _k
(m_- ^0(x ,k)  m_+ ^0(y ,k))    \big | \le   \frac{C }{   |k-\pi
\ell _{m }|
 +|g_m| }.
 $

\endproclaim
  By Lemmas 4.1-3 and Lemma 3.1 we conclude:

\proclaim{Lemma  4.4} $\exists $   $ \Bbb D\subset (0,+\infty )$,
bounded discrete, such that $\forall$ $\mu \in (0,+\infty )
\backslash \Bbb D$ and and for any $n_0 $ bands then there exists a
$C=C(\mu ,n_0)
>0$ such that for all $x,y$ and for all $t\ge 0$ we have $ \left |
\sum _{|n|\le n_0} K_n(t,x,y)\right |\le C  \langle t \rangle
^{-\frac{1}{3}} .$
\endproclaim
As a consequence of Lemma 4.4, in order to prove Proposition 2.2 it
is enough to look at  $K^n$ in (4.1) with large $n$. It is not
restrictive to sum over   $n\gg 1$. We split further in (4.1). We
introduce a smooth, even, compactly supported cutoff $\chi
_{0}(t)\in [0,1]$ with $\chi \equiv 0$ near 1 and $\chi _{0}=1$ near
$0$. Set $\chi _{1} =1-\chi _{0} $. For $c\gg 1$ fixed we split each
$K^n$ in (4.1) as $K^n=\sum _1^5K^n_j $ partitioning the identity in
$\sigma _n=[a_{n}^+, a^-_{n+1}]$,
$$\aligned &  1_{\sigma _n}(w)= \chi _0 (\frac{w - a_{n}^+ }{c|g_n|}  ) +  \chi _1 (\frac{w - a_{n}^+ }{c|g_n|}  )
  \chi _0(\frac{w - a_{n}^+ }{
 |g_n|^{ \frac{1}{4}}} )  \\&   +  \chi _1 (\frac{w
- a_{n}^+ }{  |g_n|^{ \frac{1}{4}}}  )   \chi _1 (\frac{ a_{n+1}^--w
}{ |g _{n+1}|^{ \frac{3}{5}}}  ) + \chi _0(\frac{ a_{n+1}^--w }{ |g
_{n+1}|^{ \frac{3}{5}}} )
  \chi _1  (\frac{
a_{n+1}^--w }{c |g _{n+1}|} ) +\chi _0(\frac{ a_{n+1}^--w }{ c|g
_{n+1}|} ) .
\endaligned
$$
  By $ c\gg 1$ we have $  \dot w
\approx 1$ for $w\in [a_{n}^++c|g_n|, a^-_{n+1}-c|g_ {n+1}|]$ ,
Lemma 3.6.

\proclaim{Lemma 4.5} $\exists$ a fixed $C>0 $ s.t. $ |K^n_1| \le C
t^{-\frac{1}{2}}| g_{n}|^{\frac{1}{2} }$ and $ |K^n_5|  \le C
t^{-\frac{1}{2}}| g_{n+1}|^{\frac{1}{2} }$.
\endproclaim

\proclaim{Lemma 4.6 } There are an $\epsilon >0$ and $C_\epsilon $
such that $ |K^n_2|  \le C_\epsilon t^{-\frac{1}{3}}| g_{n}|^{
\epsilon}$ .
\endproclaim

\proclaim{Lemma 4.7} There are an $\epsilon >0$ and $C_\epsilon $
such that $ |K^n_4|  \le C_\epsilon t^{-\frac{1}{3}}  | g_{n+1}|^{
\epsilon}  $.
\endproclaim

\noindent Lemmas 4.5-7 imply $ \sum _n \sum _{j \neq 3,j =1}^5
|K^n_j (t,x,y)| \le C \max \{ t^{-\frac{1}{3}}, t^{-\frac{1}{2}}\}
.$ Turning to $K^n_3 $,  set $K_3=\sum _n K _3^n $. The following
lemma completes the proof of Proposition 2.2:

\proclaim{Lemma  4.8} There is a fixed $C$ such that $|   K _3
(t,x,y)| \le C\langle t \rangle ^{-\frac{1}{3}}   .$
  \endproclaim
We   prove Lemmas 4.5-7 in \S 5 and  Lemma 4.9 in \S 6.

  \head \S 5 Proof of   Lemmas 4.5-7 \endhead

For all the $j\neq 3$ and for $\psi (k)$ the corresponding cutoff,
we consider $$H^n_j(t,x,y)=
  \int   _{[\ell _n \pi ,\ell _{n+1}\pi ]}e^{ -i(x-y)k \pm it\eta (k) }
    m _- ^0(x,k)  m  _+ ^0(y,k) \eta ^{- \frac{3}{2}}(k) \psi (k)dk. $$
Lemma 4.5 is an immediate consequence of:

 \proclaim{Lemma 5.1} $\exists$ a  fixed $C $ such
that $ |H^n_1|  \le C  t^{-\frac{1}{2}}| g_{n}|^{\frac{1}{2} }$ and
$ |H^n_5|  \le C  t^{-\frac{1}{2}}| g_{n+1}|^{\frac{1}{2} }$.
\endproclaim
{\it Proof.} We will prove the $j=1$ case. Recall from formula (3.1)

$$\ddot \eta  =  2^{-1} {\ddot E} (E+\mu
)^{-\frac{1}{2}}     - 4^{-1}{\dot E^2}{ (E+\mu )^{-\frac{3}{2}}}.$$
For $0\le w -a^{+}_{n}\lesssim |g_n|$    by Lemma 3.6 we have $0\le
\dot w\lesssim (w
-a^{+}_{n})^{\frac{1}{2}}|g_n|^{-\frac{1}{2}}\lesssim 1$ and so in
particular $\dot E^2\lesssim \ell _n^2$. By Lemma 3.5 we have
$|\ddot E|\ge c \ell _n |g_n|^{-1}$ for some fixed $c>0$. Hence
$|\ddot \eta | \gtrsim |g_n|^{-1}$. Then, by Lemmas 4.1-3 we obtain

$$|H^n_1(t,x,y)| \le \frac{C\sqrt{| g_{n}|}}{ \sqrt{  t}} \int
_{a^{+}_{n}} ^{a^{+}_{n}+c|g_n|} \frac{  \frac{dk}{dw} dw}{ |k-\pi
\ell _{n}|
 +|g_n| }  .
$$
We have $\frac{dk}{dw}\approx  {\sqrt{|
g_{n}|}}{(w-a^{+}_{n})^{-\frac{1}{2}}}$ and so $|k-\pi \ell
_{n}|\approx \sqrt{| g_{n}|} \sqrt{w-a^{+}_{n} }.$ Hence

$$|H^n_1(t,x,y)| \le \frac{C  _1 }{ \sqrt{  |g_n|t}} \int
_{a^{+}_{n}} ^{a^{+}_{n}+c|g_n|}  \frac{\sqrt{|
g_{n}|}}{\sqrt{w-a^{+}_{n}}}dw \le \frac{C_ 2 \sqrt{|g_n|}  }{
\sqrt{ t}}.
$$
The argument for  $H^n_5$ is the same, by $\dot E^2\lesssim \ell
_{n+1}^2$ and $|\ddot E|\ge c \ell _{n+1} |g_{n+1}|^{-1}$.

\bigskip

Lemma 4.6 is an immediate consequence of the following lemma:

 \proclaim{Lemma 5.2} There is $C> 0$ such that $ |H^n_2| \le C
\min\{ \langle \ell _n\rangle ^{ \frac{3}{2}}t^{-\frac{1}{2}}\log
(1/|g_n|),  | g_{n}|^{\frac{1}{4}} \}$.
\endproclaim
{\it Proof.} $H^n_2$ is defined by an integral for $w\in
[a^{+}_{n}+c |g_n|, a^+_{n }+ |g_{n }|^{\frac{1}{4}}]$. We claim we
have $  |\ddot  \eta | \gtrsim \langle k\rangle ^{ -3} . $ Assume
this inequality. By Lemma 3.6 we have $\frac{dk}{dw}\approx 1$,
$w-a^+_n \approx k-\pi \ell _{n } $.  So by Lemmas 4.1-3, by $ k-\pi
\ell _{n} \gtrsim  |g_n|$ and proceeding as in  Lemma 5.1
$$ \aligned & |H^n_3|\le C _1t^{-\frac{1}{2}}   \langle \ell _n\rangle ^{ \frac{3}{2}}
\int _ {a^+_{n }+ c|g_{n }| } ^{ a^+_{n }+  |g_{n }|^{\frac{1}{4}}}
 \frac{  \frac{dk}{dw} dw}{ |k-\pi
\ell _{n}|
 +|g_n| }
    \le \\& \le C_{2} t^{-\frac{1}{2}}   \langle \ell _n\rangle ^{ \frac{3}{2}}  \int _ {a^+_{n }+
c|g_{n }| } ^{ a^+_{n }+  |g_{n }|^{\frac{1}{4}}} \frac{   dw}{
w-a^+_{n }
  }  \le C_3 t^{-\frac{1}{2}}   \langle \ell _n\rangle ^{ \frac{3}{2}}
 \log \frac 1{|g_n|}  . \endaligned
$$
By Lemma 4.2 we have also $|K^n_3|\le C |g_n|^{\frac{1}{4}}.$ To
prove $ |\ddot  \eta | \gtrsim \langle k\rangle ^{ -3}   $  we write
$ 4\ddot \eta  =( 2(E+\mu )\ddot E-\dot E^2) { (E+\mu
)^{-\frac{3}{2}}} $ with $E= k^2+2Q_0+O(k^{-2})$, $\dot
E^2=4k^2+O(1/k)$:
$$\ddot \eta = \frac{(2\ddot E-4)k^2+2\ddot E(2Q_0+\mu ) +O(1/k)}{4(E+\mu )^{ \frac{3}{2}}}
.$$ For $w\in [a^{+}_{n}+c |g_n|, a^+_{n }+ |g_{n }|^{\frac{1}{4}}]$
we have $k''<0$ and so as in Lemma 3.3
$$ \ddot E=2(\dot w)^2-2w (\dot w)^3k''\ge 2(\dot
w)^2=2(1-2Q_0k^{-2}+O(k^{-3})).$$ So we get $\ddot \eta \ge (  {\mu
+O(k^{-1})})  {(E+\mu )^{-\frac{3}{2}}}  $ and our claim is proved.

Lemma 4.7 is an immediate consequence of the following lemma.
\proclaim{Lemma 5.3} There is a $C  $ s.t. $ |H^{n-1}_4| \le C
t^{-\frac{1}{3}} | g_{n-1}|^{ 1/30}$.
\endproclaim
{\it Proof.} $H^{n-1}_4$ is defined by an integral for $w\in [
a^-_{n  }- |g_{n }|^{\frac{3}{5}},
         a^{-}_{n }- c|g_n|]$.
By Lemma 4.3 \cite{C}  we have $  |\dddot  \eta | \gtrsim |g_{n
}|^{-\frac{1}{10}}. $   By Lemma 3.6 we have $\frac{dk}{dw}\approx
1$, $ a^-_n -w\approx  \pi \ell _{n } -k$, and so by Lemmas 4.1-3
and proceeding as in  Lemma  5.1
$$ |H^n_4|\le C t^{-\frac{1}{3}}    |g_{n }|^{ \frac{1}{30}}  \int ^
{a^-_{n  }- \ell _n |g_{n }| } _{ a^-_{n  }-  |g_{n
}|^{\frac{3}{5}}} \frac{ dw}{  a^-_{n }-w
  }  \le C_1 t^{-\frac{1}{3}}  |g_{n }|^{ \frac{1}{30}}
 \log \frac 1{|g_n|}  .
$$

\head \S 6 Proof of Lemma 4.8\endhead

\noindent  For $K_3=\sum _n K _3^n $ we  show $| K _3 (t,x,y)| \le
C\langle t \rangle ^{-\frac{1}{3}}    $ for $C$ fixed  by reducing
to the flat case of Lemma 5 \cite{MSW}, whose proof permeates this
section. Set $ \chi _{int} (k)=\sum _n \chi _1 (\frac{w - a_{n}^+ }{
|g_n|^{ \frac{1}{4}}}  ) \chi _1 (\frac{ a_{n+1}^--w }{ |g _{n+1}|^{
\frac{3}{5}}}  ) $ supported inside the union of sets $ a_{n}^+ +
|g_n|^{ \frac{1}{4}} \le w \le a_{n+1}^-- |g _{n+1}|^{ \frac{3}{5}}
$. Then for $R=x-y$

 $$K_3(t,x,y)=2
  \int   _{0}^\infty \cos (Rk  ) \sin (t\eta (k)  )  \eta ^{-\frac{3}{2}}(k)
  {m _- ^0(x,k)}{ m  _+ ^0(y,k)} \chi _{int} (k) dk .$$
  We split the integral between $[0,t]$ and $[t,\infty )$. By Lemma
  4.2 and by $\eta (k) \approx \langle k \rangle $  the $[t,\infty )$ integral has absolute
   value less than  $C\langle t \rangle ^{-\frac{1}{2}} $ for a fixed $C>0$.
Next write $$\aligned & h_\pm (k)= t   \eta (k)\pm Rk \\& I_\pm
(t)=\int _{0}^t e^{i h_\pm (k)}     \eta ^{-\frac{3}{2}}(k)
  {m _- ^0(x,k)}{ m  _+ ^0(y,k)} \chi _{int} (k) dk \endaligned  \tag 6.1$$
It is not restrictive to assume $R=x-y>0$. We start with $I_-(t)$.

\proclaim{Lemma  6.1 } There is a fixed $C$ such that $|  I_- (t)|
\le C\langle t \rangle ^{-\frac{1}{3}}  .$
  \endproclaim
{\it Proof.}   The proof ends in  Lemmas 6.16.  We set $ I_-(t)=
I_1(t)+I_2(t) $ with
$$\aligned &
I_{1}( t)=\int   _{0}^t e^{i h(k)}     \eta ^{-\frac{3}{2}}(k)
  \left ( m _- ^0(x,k)  m  _+ ^0(y,k)-1 \right ) \chi _{int} (k) dk
\\& I_{2}( t)=\int   _{0}^t e^{i h(k)}     \eta ^{-\frac{3}{2}}(k)
    \chi _{int} (k) dk
\endaligned \tag 6.2
$$

To prove Lemma 6.1 we   use: \proclaim{Lemma  6.2} In   $supp (\chi
_{int})\cap [0,t]$ we   have for fixed constants:

{\item {(1)}} $ 0<1-\dot \eta (k)<Q_0 k^{-2}+O(k^{-3}) $;

{\item {(2)}} $ \ddot \eta (k)\gtrsim   \langle k \rangle ^{-3} $
and so $|\ddot h (k)|= t |\ddot \eta (k)|\gtrsim t \mu \eta
^{-3}\approx t \langle k \rangle ^{-3}$.

 {\item {(3)}} There is a fixed   $C>0$   such that
for any $n$ sufficiently large we have
$$|\dot E(k) - 2k| \le \frac{C}{\langle k \rangle
  ^2} \quad \text{for} \quad   a^+_n + \ell ^3_n|g_n| \le w   \le a^-_{n+1}-
  \ell ^3_{n+1}|g_{n+1}| . $$

  {\item {(4)}} There are fixed constants $C>0,$ $C_1>0 $ and
$c>0$ such that for  any $n$  and any $k\in [\ell _n\pi , \ell
_{n+1} \pi ] $, that is for any $w\in [a^+_n , a^-_{n+1}]  $, we
have:

$$\aligned    &     a^+_n +c|g_n| \le w\le     a^-_{n+1}-C _1\ell
_{n+1} ^{\frac{1}{3}}|g_{n+1 }|^{\frac{2}{3}}   \Rightarrow \ddot{E}
\approx \frac 12 +\frac{\ell _n|g_n|^2}{ |w-a_n^+|^{3}} .
 \endaligned $$

\endproclaim
{\it Proof.} For (3)  see Lemma 7.1 \cite{C}, for (4) see Lemmas 4.2
and 7.4 \cite{C}. (2) is proved as in Lemma 5.2.   By construction
$supp (\chi _{int})\cap [0,\infty ) \subset [\widetilde{k} ,\infty
)$ for some $\widetilde{k}\approx n_0\gg 1$. The following for $k\gg
1$, which uses (3) here and (4) Lemma 3.3, proves (1):

$$\dot \eta  = \frac{\dot E}{2(E+\mu )^{\frac{1}{2}}}= \frac{2k+O(k^{-2})}{ 2(k^2+2Q_0+O(k^{-2}))^{\frac{1}{2}}}
=1-Q_0 k^{-2}+O(k^{-3}).$$

\bigskip
The following lemma coincides with Lemma 4.7 \cite{C}, with the
proof scattered in Lemmas 5.2, 7.1, 7.4 and 7.6 \cite{C}:
\proclaim{Lemma 6.3} In the support of $\chi _{int} $ we have
$w\approx k$, $\dot w= 1+O(k^{-2})$, $\ddot w= O(k^{-3})$. We can
extend $w$ from the support of $\chi _{int} $ to the whole of $\Bbb
R$ so that the extension (which we denote again with $w$) satisfies
the same relations and is an odd function.
\endproclaim
Thanks to Lemma 6.3  we obtain:

\proclaim{Lemma 6.4} We can extend $\eta (k)=\sqrt{w^2(k)+\mu}$   to
all $\Bbb R$ so that there are fixed positive $ c_1 $, $c_2$ so that
$ \ddot \eta (k) \ge c_1 \langle k \rangle ^{-3} $ and $ |1-\dot
\eta (k)|\le c_2 \langle k \rangle ^{-2} $,  and   positive $ c_3 $,
$c_4$ such that in $\Bbb R \backslash [-1,1]$,  $ c_3\ge \dot \eta
(k) \ge
 c_4$. Furthermore, from $\dot \eta = \frac{w\dot
 w}{\sqrt{w^2+\mu}}$ where $\mu \ge \mu _0>0$, from Lemma 6.3,   $n_0=n_0(\mu _0)$ can be chosen and the extension in
 Lemma 6.3 be done so that $|\dot \eta |<1$ in $\Bbb R$.
 \endproclaim
 In the rest of the  paper by $\eta (k)$ we
will mean this extension  and we will set $h(k)=t\eta (k)-Rk$. We
have:

\proclaim{Lemma  6.5} Consider the $h(k)$ just introduced. {\item
{(1)}}If $t\le R$ then $| \dot h (k) |=R- t
 \dot \eta (k)\ge c t |k| ^{-2}$ for a fixed $c>0$.   {\item
{(2)}} If $t> R$ then   $\dot h(k)$ has
 exactly one  zero in $[0,+\infty )$ which we denote by $k_0$. {\item
{(3)}} In case (2), if $k_0>2$,   for  $1\le k<k_0/2$ and for
$k>2k_0$ we have $| \dot h (k) |  \ge c t |k| ^{-2}$. If $k_0\le 2$
for $k>2k_0$ we have $| \dot h (k) |  \ge c t |k| ^{-2}$.
\endproclaim
{\it Proof.}  $|\dot \eta |<1$ implies (1). Consider $t> R$. By
$\dot h=  {tw\dot
 w}(w^2+\mu) ^{-\frac{1}{2}}-R$ we have $\dot h (0) = -R$  and by Lemma 6.3
$\dot h\approx t-R>0$ for $k\to \infty$. So there is a zero which by
  $\ddot h= t\ddot \eta \ge c_1 t \langle k \rangle
^{-3}>0 $   is unique. We denote it by $k_0$. This gives us (2).
 We set $[a,b]=[1,t]\cap [k_0/2,2k_0]$. For $k\in (1,a)$
 by Lemma 6.5
    we have
 $ \dot h (k)< \dot h (2k)<0$.
 So for some $\widetilde{k}\in [k ,2k]$
$$|\dot h (k)|>\dot h (2k)-\dot h (k)=\ddot h ( \widetilde{k})k>    ct   k^{-2}.$$
For $k>b\ge 2k_0$ we have $\dot h ( k)>\dot h (k/2)>0$ and for some
$\widetilde{k}\in [k/2, k]$
 $$|\dot h (k)|>\dot h ( k)-\dot h (k/2)=\ddot h ( \widetilde{k})k/2>   ct   k^{-2}.$$
Lemmas 4.1 and 6.2-4 imply:

\proclaim{Lemma  6.6} Let $\dot H (k)=e^{i h(k)}$ with $H(0)=0$.
Then for a fixed $c>0$    we have $|H(k)|\le c
t^{-\frac{1}{2}}\langle k \rangle ^{\frac{3}{2}} $ for all $k\in
[0,t]$.
\endproclaim
Next, we have the following analogue of Lemma 5 \cite{MSW}:
\proclaim{Lemma  6.7} For $|g(k)|=O(\langle k \rangle
^{-\frac{5}{2}})$ we have $$ \left | \int _0^t H(k) g(k) dk \right
|\le C \langle t \rangle ^{-\frac{1}{2}}  .\tag 1$$ \endproclaim

{\it Proof.}   By Lemma  6.6, $|H(k)|\lesssim \langle t \rangle
^{-\frac{1}{2}} $ for $|k|\le 2$. If $R\ge t$ by Lemma 6.4 we have
$|\dot h(k)|\ge c t k^{-2}$. Then by Lemma 4.1 for $k\ge 1$ we have
$|H(k)-H(2)|\le c t^{-1} k^{2} .$ By $|H(2)|\le C \langle t \rangle
^{-\frac{1}{2}} $ and by $|g(k)|=O(\langle k \rangle
^{-\frac{5}{2}})$ we obtain (1). If $R< t$
 by Lemma 6.4  $\dot h(k)$ has   one zero  which we denote by $k_0$.
 If $ k_0\le 2$ we can repeat the above argument. If $ k_0>  2$
 set $[a,b]=[1,t]\cap [k_0/2,2k_0]$ as in Lemma 6.5.   For $k\in (1,a)$ and $p<k$
  by Lemma 6.5
$$|\dot h (p)|>  {ct}{ p ^{-2}}> {Ct}{ k^{ -2}}.$$
By Lemmas 4.1 and 6.5 $|H(k)|\le ck^2t^{-1}$. Then we get
$$\left | \int _1^a H(k) g(k) dk \right |\le C  t   ^{-\frac{1}{2}}
. $$ For $k>p>b\ge 2k_0$ by Lemma 6.5
 $$|\dot h (p)|>  {ct}{ p ^{-2}}> {Ct}{ k^{ -2}}.$$
By Lemmas 4.1 and 6.6 $|H(k)-H(b)|\le ck^2t^{-1}$. Then we get

$$\left | \int _b^t \left (  H(k)-H(b)\right )  g(k) dk \right |\le C  t   ^{-\frac{1}{2}}
. $$ By Lemma 6.6

$$\left | \int _b^t  H(b)   g(k) dk \right |\le C  t
^{-\frac{1}{2}} \langle b \rangle  ^{\frac{3}{2}} b^{-\frac{3}{2}} .
$$ Finally

$$\left | \int _a^b H(k) g(k) dk \right |\le C   t   ^{-\frac{1}{2}}
\int _{\frac{k_0}{2}}^{2k_0} \frac{dk}{k}\le C  t
^{-\frac{1}{2}}2\log 2   .$$

\proclaim{Lemma  6.8} There is a fixed $C$ such that for the
$I_2(t)$   in (6.2),  $| I_2(t)| \le C\langle t \rangle
^{-\frac{1}{3}}  .$
  \endproclaim
{\it Proof.}  Write $\chi _{int} (k)=1-\chi _{ext} (k)$ and
correspondingly $I_{2}( t)=I_{21}( t)-I_{22}( t) $ with
 $I_{21}( t)=\int   _{0}^t e^{i h(k)}     \eta ^{-\frac{3}{2}}(k)
      dk $  and $I_{22}( t)=\int   _{0}^t e^{i h(k)}   \chi _{ext} (k)  \eta ^{-\frac{3}{2}}(k)
      dk $.
Then:

\proclaim{Lemma  6.9} There is a fixed $C$ such that $|  I_{21}( t)|
\le C\langle t \rangle ^{-\frac{1}{2}}  .$
  \endproclaim

\proclaim{Lemma  6.10} There is a fixed $C$ such that $|  I_{22}(
t)| \le C\langle t \rangle ^{-\frac{1}{3}}  .$
  \endproclaim
{\it Proof of Lemma 6.9.} We have
$$I_{21}( t)= \frac{3}{2}\int   _{0}^t  H(k)      \eta ^{-\frac{5}{2}}(k)
\dot \eta (k)      dk+H(t)  \eta ^{-\frac{3}{2}}(t)  .$$ $|H(t)|\le
t$ and $\eta (t)\approx \langle t \rangle  $ imply $|H(t) \eta
^{-\frac{3}{2}}(t)|\lesssim \langle t \rangle ^{-\frac{1}{2}}$.
  We have $\eta (k)\approx k$  and   $ |\dot \eta |\lesssim 1 $.
So $g(k):=\eta ^{-\frac{5}{2}}(k) \dot \eta (k) =O(\langle k \rangle
^{-\frac{5}{2}}) $. Then Lemma 6.7 implies Lemma 6.9.

  {\it Proof of Lemma 6.10.} For $J_n= [\ell _n\pi -
C|g_n|^{\frac{3}{5}},
    \ell _n\pi +   C |g_n|^{\frac{1}{4}} ]
   $   for some fixed $C \gtrsim 1$, set
$$ \aligned & I_{22}( t)=\int   _{0}^t e^{i h(k)}     \eta ^{-\frac{3}{2}}(k)
    \chi _{ext} (k)  dk= -I_{221}( t)-I_{222}( t)\, , \quad
I_{221}( t):=   \sum _{n=0}^\infty
    \int _{[0,t]\cap J_n} \\& \times H_n(k) \chi _{ext} (k) \frac{d}{dk}  \eta ^{-\frac{3}{2}}(k)
      dk\, , \quad I_{222}( t):= \sum _{n=0}^\infty
    \int _{[0,t]\cap J_n} H_n(k)  \eta ^{-\frac{3}{2}}(k) \dot \chi _{ext} (k)
      dk,\endaligned
$$ with $\dot H _{n}(k)=e^{i h(k)}$, $H_{n}(\ell _n)=0$.  By $\ddot h (k)= t \ddot \eta  (k)$,
   $c_1   \langle k \rangle ^{-3}\le
\ddot \eta (k)$, $|\dot \eta |<1$ and Lemma 6.6 which implies
$|H_n(k)|\le Ct^{-\frac{1}{2}}
 \langle k \rangle ^{ \frac{3}{2}} $ for $C$ fixed,
$$|H_n(k) \chi _{ext}  (k)  \eta ^{-\frac{5}{2}}  (k)\dot \eta (k)|\le C  t^{-\frac{1}{2}}
 \langle k \rangle ^{ \frac{3}{2}} \langle k \rangle ^{ -\frac{5}{2}} $$
so $|I_{221}( t)|\le C t^{-\frac{1}{2}} \sum _{n=1}^ {[t]} \langle n
\rangle ^{ -1}  \lesssim   t^{-\frac{1}{2}}|\log t|.$ By
$|H_n(k)|\le C \min \{ t^{-\frac{1}{2}}
 \langle k \rangle ^{ \frac{3}{2}}), |g_n|^{\frac{1}{4}}\}$

$$ \int _{[0,t]\cap J_n}| H_n(k)  \eta ^{-\frac{3}{2}}(k) \dot \chi _{ext}
(k)|
      dk\le  \min \{ t^{-\frac{1}{2}}
 , |g_n|^{\frac{1}{4}} \langle \ell _n \rangle ^{ -\frac{3}{2}} \} .
$$
This by Theorem 2.3 implies $|I_{222}( t)|\le C t^{-\frac{1}{3}}.$

\proclaim{Lemma  6.11} There is a fixed $C$ such that $|  I_1(t)|
\le C\langle t \rangle ^{-\frac{1}{3}}  .$
  \endproclaim
{\it Proof.} We fix  some small $\varepsilon >0$ and   split $I_{1}(
t)=I_{11}( t)+I_{12}( t) $ with

$$\aligned & \widetilde{\chi } _{int}(k):=\sum _n \chi _{1}\left ( \varepsilon^{-1} (k - \pi
\ell _{n} ) \right ) \chi _{1} \left (\varepsilon^{-1} (\ell _{n+1}
\pi -k ) \right )\\& I_{11}( t):= \int   _{0}^t e^{i h(k)} \eta
^{-\frac{3}{2}}(k)
      \left ( m _- ^0(x,k)  m  _+ ^0(y,k)-1 \right )   \widetilde{\chi } _{int}(k) dk
      \\& I_{12}( t):= \int   _{0}^t e^{i h(k)}     \eta ^{-\frac{3}{2}}(k)
      \left ( m _- ^0(x,k)  m  _+ ^0(y,k)-1 \right ) \chi _{int} (k) (1-\widetilde{\chi } _{int}(k)) dk
      .\endaligned
$$
\proclaim{Lemma 6.12} There is a fixed $C$ such that $|I_{11} (t)
|\le C \langle t\rangle ^{-\frac{1}{2}}$.
 \endproclaim
{\it Proof.} We   have  $ m _- ^0(x,k)  m _+ ^0(y,k)-1= O( k^{-1})$
and $ \partial _k \left ( m _- ^0(x,k) m _+ ^0(y,k)\right ) = O(
k^{-1}) $ in the support of $\widetilde{\chi } _{int} (k).$ Then
Lemma 6.7 implies Lemma 6.12.

\bigskip
\proclaim{Lemma 6.13} There is a fixed $C$ such that $|I_{12} (t)
|\le C \langle t\rangle ^{-\frac{1}{3}}$.
 \endproclaim
{\it Proof.} We consider $I _{12} (t )=\sum _n I _{12} ^{n}(t )$, $
I_{12}^{n} (t )
   :=\int _{\ell _{n}\pi}^{\ell _{n+1} \pi }
e^{i  h ( k)} f(k) dk$
$$\aligned &  f(k):=\Psi _n(k)  \eta ^{-\frac{3}{2}} (k) \left ( m _-
^0(x,k) m _+ ^0(y,k)-1 \right ) \text{ where } \\& \Psi _n(k):= \chi
_{1} (\frac{w - a_{n}^+ }{ |g_n|^{ \frac{3}{5}}} ) \chi _{1} (\frac{
a_{n+1}^--w }{ |g _{n+1}|^{ \frac{1}{4}}} )\chi _{0}   (\frac{k -
\pi \ell _{n } }{\varepsilon }  ) \chi _{0}  (\frac{\ell _{n+1} \pi
-k }{\varepsilon}  )  .\endaligned
$$
  Observe
that $ \Psi _n(k)=  \Psi _{n 1}(k) +\Psi _{n 2}(k)$  with $\Psi _{n
1}(k)$ supported in $ |g_n|^{ \frac{1}{4}}\lesssim k-\pi \ell _{n
}\lesssim \varepsilon$ and with $\Psi _{n 2}(k)$ supported in $
\varepsilon \gtrsim    \pi \ell _{n+1}-k\gtrsim  |g _{n+1}|^{
\frac{3}{5}} .$ Correspondingly write $f=f_1+f_2$ and $I _{12}^{n}=
I _{12}^{n1}+
   I_{12}^{n2} .$ We have:
   \proclaim{Lemma 6.14} For a fixed $C$ and  for  $j=1,2$: $\sum
   _n|I
   _{12 }^{nj}(t )|\le C   \langle t  \rangle  ^{-\frac{1}{2}} |\log   t
    |  .$
   \endproclaim
{\it Proof.} We  focus on $I _{12}^{n1}$, the proof for  $I
_{12}^{n2}$ being almost the same. We have
$$ I _{12}^{n1} (t,x,y)
   =\int _{\ell _{n }\pi}^{(\ell _{n }+\frac{1}{2}) \pi } \dot H  _n(k) f_1(k)dk
    \text{ with } H_n(k)= \int _{\ell _{n }\pi}^{k}e^{i  h ( k^\prime
   )} dk^\prime .$$
   For $I
_{12}^{n2}$ the proof is the same but with  $H_n(k)= \int _{\ell
_{n+1} \pi}^{k}e^{it h ( k^\prime
   )} dk^\prime $.  We get
   $$\aligned & I _{12}^{n1} (t )
   =-I _{121}^{n1} (t )-I _{122}^{n1} (t ) \text{ with }  \\& I_{121}^{n1} (t
   )
 =\int _{\ell _{n}\pi}^{(\ell _{n}+\frac{1}{2}) \pi
} H_n(k)  \partial _k\left ( \Psi  _{n1} (k) \eta ^{-\frac{3}{2}}(k)
\right )
  \left ( m _- ^0(x,k) m _+ ^0(y,k)-1 \right ) dk\\& I _{122}^{n1}
(t )=\int _{\ell _{n}\pi}^{(\ell _{n}+\frac{1}{2}) \pi } H_n(k)\Psi
_{n1}   (k)  \eta ^{-\frac{3}{2}}(k)   \partial _k \left ( m _-
^0(x,k) m _+ ^0(y,k)-1 \right )   dk .\endaligned
 $$

  \proclaim{Lemma 6.15}   There is a fixed $C>0$ such that $ \sum _{n=1}^{\infty }|I _{121}^{n1} (t
)|\le C \langle t \rangle ^{-\frac{1}{3}} .$
\endproclaim
{\it Proof.}
 Set
$$\partial _k\left ( \Psi  _{n1} (k) \eta ^{-\frac{3}{2}}(k)
\right )=  -\frac{3}{2} \Psi  _{n1} (k)\eta ^{-\frac{5}{2}}(k) \dot
\eta
 (k)+ \dot \Psi _{n1} (k)\eta ^{-\frac{3}{2}}(k)=:a(k)+b(k).$$ We
 have $a(k)=O(k^{-\frac{5}{2}})$ and
$$|b(k)|\le C ( \varepsilon ^{-1}k^{-\frac{3}{2}} \chi _{I_1}(k) +
  |g_n|^{
-\frac{ 1}{4}}k^{-\frac{3}{2}} \chi _{I_2}(k) )$$ with length of
$I_1\approx \varepsilon $ and length of $I_2\approx  |g_n|^{  \frac{
1}{4}} $. Recall by Lemma 4.2, $$ \left ( m _- ^0(x,k) m _+
^0(y,k)-1 \right ) =O(k^{-1}).$$ For $t\le R$ then $\dot h \ge ct
\langle k\rangle ^{-2}$. From the estimates on $a(k)$ and $b(k)$ we
get $|I _{121}^{n1} (t )|\le C t^{-1} \ell _n^{-\frac{1}{2}}.$
Summing up over $\ell _n \lesssim t$ we get much less than $t
^{-\frac{1}{3}}$. For $t> R$ and for the critical point $k_0>2$
(otherwise proceed as above) distinguish between  two cases

Case 1: $\pi \ell _n$ outside $[k_0/2-\pi, 2k_0+\pi ]$. Then   $
|\dot h| \ge ct \langle k\rangle ^{-2}$ by Lemma 6.5 and $|I
_{121}^{n1} (t )|\le C t^{-1} \ell _n^{-\frac{1}{2}} $. Summing up
over $\ell _n \le t$ we get $O(t^{-\frac{1}{2}})$.

Case 2:  $\pi \ell _n$ inside $[k_0/2-\pi, 2k_0+\pi ]$. Then by
Lemma 6.4 we have $|H_n(k)|\le C t^{-\frac{1}{2}}\langle k \rangle
^{\frac{3}{2}} $ and $|I _{121}^{n1} (t )|\le C t^{-\frac{1}{2}}
\ell _n^{-1} $. Summing up over $\ell _n \approx k_0$ we get
$O(t^{-\frac{1}{2}})$.

\proclaim{Lemma 6.16} There is a fixed $C>0$ such that $ \sum
_{n=1}^{\infty }|I _{122}^{n1} (t )|\le C \langle t \rangle
^{-\frac{1}{3}} .$
\endproclaim
{\item {(1)}} Suppose $t> R$. Then $|\dot h(k)|\ge  ct \langle k
\rangle ^{-2}$. For a fixed $C$ by Lemma 4.1
$$\big | \int
_{\ell _{n}\pi}^{k}e^{i   h ( k^\prime
   )} dk^\prime \big |  \le \min \{ C   \langle k \rangle ^{ 2}
t^{ -1} ,   |k- \pi \ell _{n}|  \} .\tag 6.3$$ Next we split $ I
_{122}^{n} (t )=\int _{\ell _n\pi}^{\ell _n \pi +  t  ^{ -1} } \dots
+\int ^{(\ell _n+\frac{1}{2})\pi}_{\ell _n \pi +  t  ^{ -1}  } \dots
$. By Lemma 4.3
$$ \aligned & \big |\int _{\ell _n\pi}^{\ell _n \pi +  t  ^{ -1} }
H_n(k)\Psi _{n1}   (k)  \eta ^{-\frac{3}{2}}(k)   \partial _k \left
( m _- ^0(x,k) m _+ ^0(y,k)  \right ) \big | \\&\le C \int _{\ell
_n\pi}^{\ell _n \pi +  t  ^{ -1} } |k- \pi \ell _n| |k- \pi \ell _n|
^{-1} \langle \ell _ n \rangle ^{ - \frac{5}{2}} dk=  C t  ^{
-1}\langle \ell _ n \rangle ^{ - \frac{5}{2}}
 \endaligned \tag 6.4
$$ and
$$ \big |\int ^{(\ell _n+\frac{1}{2})\pi}_{\ell _n \pi +  t  ^{ -1} } \dots  \big |
 \le   \int ^{(\ell _n+\frac{1}{2})\pi}_{\ell _n \pi +  t  ^{ -1} }  C t  ^{
-1}\langle k \rangle ^{ 2}|k- \pi \ell _n| ^{-1} \langle  \ell _n
\rangle ^{ - \frac{5}{2}}\tag 6.5
$$
and so $ |I _{122}^{n} (t)| \le C\langle  \ell _n \rangle ^{ -
\frac{1}{2}}t ^{ -1} \log t .$  Then $\sum |I _{122}^{n} (t)| \le C
t ^{ -\frac{1}{2}} \log   t.$

{\item {(2)}} Suppose $t<R$. Then there is  a unique $k_0> 0$ with
$\dot h(k_0)=0$. If $k_0\le 2$ we have $ \dot h(k) \ge  ct \langle k
\rangle ^{-2}$ in the support of the integrands and we can apply the
argument in (1). If $k_0>2$  set $[a,b]=[1,t]\cap [k_0/2,2k_0]$.
Then consider $\ell _n \le a-\pi /2$.  Then for    a fixed $C$ we
get (6.3) and by proceeding as in the case $t>R$ we can split again
and obtain estimates (6.4-5). Same is true for $\ell_n\ge b$.
Summing up over all these $\ell _n\lesssim t$ we get $\sum |I
_{122}^{n} (t)| \le C t ^{ -\frac{1}{2}} \log   t.$ For $a-\pi /2<
\ell _n< b$ for a fixed $C$
$$\big | \int
_{\ell _{n}\pi}^{k}e^{it  h ( k^\prime
   )} dk^\prime \big |  \le \min \{ C t
^{ -\frac{1}{2}} k ^{\frac{3}{2}},   |k- \pi \ell _{n}|  \} .$$ Next
we split
$$I _{122}^{n} (t )=\int _{\ell _n\pi}^{\ell _n \pi +  t^{-\frac{1}{2}} \ell _{n} ^{\frac{3}{2}} }
\dots +\int ^{(\ell _n+\frac{1}{2})\pi}_{\ell _n \pi +
t^{-\frac{1}{2}} \ell _{n} ^{\frac{3}{2}}  } \dots
$$
But now
$$ \big |\int _{\ell _n\pi}^{\ell _n \pi +  t^{-\frac{1}{2}} \ell _{n} ^{\frac{3}{2}} }
\dots  \big | \le C \int _{\ell _n\pi}^{\ell _n \pi +
t^{-\frac{1}{2}} \ell _{n} ^{\frac{3}{2}}  } |k- \pi \ell _n| |k-
\pi \ell _n| ^{-1} \langle \ell _ n \rangle ^{
-\frac{5}{2}}=Ct^{-\frac{1}{2}} \ell _{n} ^{-1}
$$ and
$$ \big |\int ^{(\ell _n+\frac{1}{2})\pi}_{\ell _n \pi + t^{-\frac{1}{2}} \ell _{n} ^{\frac{3}{2}} } \dots  \big |
 \le   \int ^{(\ell _n+\frac{1}{2})\pi}_{\ell _n \pi + t^{-\frac{1}{2}} \ell _{n} ^{\frac{3}{2}} }
    t^{-\frac{1}{2}}  |k- \pi \ell _n| ^{-1} \langle  \ell _n \rangle ^{ -1}
$$
and so $ |I _{122}^{n} (t)| \le C\langle  \ell _n \rangle ^{ -1}t ^{
-\frac{1}{2}} \log t .$  Then $\sum |I _{122}^{n} (t)| \le C t ^{
-\frac{1}{2}} \log   t.$

To complete the proof of Lemma 4.9 we have   to prove the following
lemma whose proof is   analogous to the proof  for $I_-(t)$ in the
easier case $t<R$ and which we skip:

\proclaim{Lemma  6.17 } There is a fixed $C$ such that $|  I_+ (t)|
\le C\langle t \rangle ^{-\frac{1}{3}}  .$
  \endproclaim

\enddocument